\documentclass[11pt]{amsart}
\usepackage{amsmath, amsthm, amssymb}
\usepackage{verbatim}
\usepackage{enumitem}
\usepackage{xcolor}
\usepackage{hyperref}
\hypersetup{colorlinks,citecolor=blue,plainpages=false,hypertexnames=false}

\newtheorem{theorem}{Theorem}[section]
\newtheorem{lemma}[theorem]{Lemma}
\newtheorem{proposition}[theorem]{Proposition}
\newtheorem{corollary}[theorem]{Corollary}
\theoremstyle{definition}

\theoremstyle{remark}
\newtheorem{remark}[theorem]{Remark}

\numberwithin{equation}{section}
\begin{document}
\title {Curvature inequalities for submanifolds with flat normal bundle}
\author{H. A. Gururaja}
\address{Department of Mathematics, Indian Institute of Science Education and Research Tirupati, Andhra Pradesh-517619. INDIA}
 \email{gururaja@iisertirupati.ac.in}


\begin{abstract}
 We study isometric immersions $f: M^n\rightarrow M^{n+k}(c), n\geq 3,$ into space forms with flat normal bundle and constant scalar curvature $R.$ Under a suitable multiplicity condition on the second fundamental form of $f,$ we prove the sharp global inequality $R> (n-1)(n-2)c$ if $c> 0$ and $R\geq n(n-1)c$ if $c\leq 0.$ We obtain a classification for submanifolds with flat normal bundle having constant scalar and mean curvatures. Finally, in the hypersurface case analogous results are proved for conformally flat hypersurfaces with constant higher order mean curvatures.
 
\end{abstract}

%%%%%%%%%%%%%%%%%%%%%%%%%%%%%%%%%%
%%%%%%%%%%%%%%%%%%%%%%%%%%%%%%%%%%%
\maketitle
\section{Introduction}
\bigskip

Let $M^{n+k}(c)$ denote the complete simply-connected Riemannian manifold of dimension $n+k$ and constant sectional curvature $c.$ Let $f: M^n \rightarrow M^{n+k}(c)$ be an isometric immersion with flat normal bundle whose second fundamental form we denote by $\alpha.$ It is well-known (e.g. \cite{DT1}, p. 35) that for each $x\in M^n$ 
there exists a positive integer $s(x)$ and unique pairwise distinct normal vectors $\eta_{1}(x), \cdots, \eta_{s(x)}(x)$ at $x,$ called the principal normal vectors of $f$ at $x,$ such that
$$T_xM^n=E_{\eta_1(x)}\oplus^{\perp}\cdots \oplus ^{\perp} E_{\eta_{s(x)}(x)},$$

\noindent where the subspaces $E_{\eta_i(x)}$ are defined by 
$$E_{\eta_i(x)}=\{X\in T_xM^n: \alpha(x) (X, Y)=\langle X, Y\rangle \eta_i (x) \ \ for\ \ all \ \  Y\in T_xM^n\}.$$

\noindent We call $f$ a proper immersion if the integer $s(x)$ is independent of $x\in M^n.$ For convenience, we say that a proper immersion $f$ has type $(l_1, l_2,\cdots, l_s)$ for some integer $s\geq 1$ if $\dim E_{n_i(x)}=l_i$ for each $1\leq i\leq s.$ \vspace{2mm}

Our first result is of global nature.\vspace{2mm}

\begin{theorem}\label{theorem}
	Let $M^n, n\geq 3,$ be a complete Riemannian manifold with constant scalar curvature $R$ and $f: M^n\rightarrow M^{n+k}(c)$ be an isometric immersion with flat normal bundle which is proper and has type $(1, n-1).$ Then the following holds.
	
	\noindent (i) $R\geq 0$ if $c=0.$ The inequality is strict if $M^n$ is compact.\\% If equality holds and $M$ is compact, then $M$ must be flat.\\
	(ii) $R> (n-1)(n-2)c$ if $c> 0.$ \\
	(iii) $R\geq n(n-1)c$ if $c< 0.$ The inequality is strict if $M^n$ is compact.
\end{theorem} 
\vspace{2mm}

 The inequalities in Theorem \ref{theorem} are optimal \cite{L}. It should be noted that the type condition $(1, n-1)$ is necessary since Okayasu \cite{O} has constructed examples of complete hypersurfaces in ${\mathbb R}^4$ having type $(2, 2)$ and constant \textit{negative} scalar curvature. We also note that if $f$ is umbilical at a point $x\in M^n$ (i.e. $f$ has precisely one principal normal vector at $x$), then it follows trivially from the Gauss equation that the scalar curvature of $M^n$ at $x$ is greater than or equal to $n(n-1)c.$ 

\begin{remark}
	In the hypersurface case, Theorem \ref{theorem} is a consequence of \cite{DD} and the study of rotation hypersurfaces of constant scalar curvature due to Leite \cite{L}.  However, our proof of Theorem \ref{theorem} is more direct and independent. Moreover, our method can be readily adopted to the setting of higher order mean curvatures of a hypersurface. 
	\end{remark}

Concerning the type condition $(1, n-1)$ we make the following remarks. In the hypersurface case, this condition is equivalent to the existence of exactly two distinct principal curvatures $\mu$ and $\lambda$ at each point of $M^n,$ with $\lambda$ having multiplicity $n-1.$ In particular, $M^n$ is conformally flat and under any condition of the form $\mu=\mu(\lambda), \lambda\neq 0,$ Do-Carmo and Dajczer \cite{DD} proved that the image $f(M^n)$ must be contained in a rotation hypersurface in $M^{n+1}(c).$ Conversely, if $M^n, n\geq 4,$ is conformally flat then by a result of Cartan \cite{C} $f$ has type $(1, n-1)$ provided it has no umbilical points. In higher codimension, Lobos and Tojeiro (\cite{LT}, proposition 2.5) showed that if $n\geq 3$ and $f: M^n\rightarrow {M}^{n+k}(c)$ has type $(1, n-1)$ with linearly independent principal normal vector fields, then $f=h\circ g,$ where $g: M^n\rightarrow M^{n+1}(c)$ is a quasiumbilical hypersurface (i.e. it has a principal curvature of multiplicity at least $n-1$ at every point) and $h: U\subset M^{n+1}(c)\rightarrow M^{n+k}(c)$ is an isometric immersion with no totally-geodesic points on an open set $U\supset g(M^n).$ In addition, on the open set of non-umbilical points of $g,$ each leaf of the umbilical foliation of $g$ is mapped into a relative nullity leaf of $h.$ Conversely, any such isometric immersion $f=h\circ g,$ with $g$ free of umbilical points, has type $(1, n-1)$ with linearly independent principal normal vector fields. However, we do not assume that the principal normals are linearly independent in this paper.\vspace{2mm}

	Suppose $M^n$ is conformally flat and $2 \leq k \leq n-3,$ and $f: M^n\rightarrow M^{n+k}(c), n\geq 5,$ is an isometric immersion that is proper and has type different from $(n).$ Although $f$ is not necessarily of type $(1, n-1)$ in this case, combining a result of Moore \cite{M} on the existence of a principal normal vector $\eta(x)$ such that $\dim{E_{\eta(x)}}\geq n-k$ for each $x\in M^n$ with a result of Dajczer, Onti and Vlachos (\cite{DOV}, theorem $1$), we infer that the possible types of $f$ are limited when the codimension $k$ is small. In particular, the following important special case in codimension two should be noted:\smallskip
	
	\textit{Let $M^n, n\geq 5,$ be conformally flat and $f: M^n\rightarrow M^{n+2}(c)$ be a proper isometric immersion that is not totally-umbilical. Then the only possible type of $f$ other than $(1, n-1)$ is $(1, 1, n-2).$}

\vspace{2mm}

%%%%%%%%%%%%%
%%%%%%%%%%%%%%%%%%%%%%%%%%%%%%%%%%%%%%%%%%%%%%%%%%%%%%%%%%%%
The main idea behind the proof of Theorem \ref{theorem} is to work with principal coordinate systems (see $\S 2$) and use a convexity argument. The same method also yields the following result which is of local character.

\begin{proposition}\label{prop}
	Let $M^n, n\geq 3,$ be a Riemannian manifold with constant scalar curvature and $f: M^n\rightarrow M^{n+k}(c)$ be an isometric immersion with flat normal bundle which is proper and has type $(1, n-1).$ Assume that the mean curvature field of $f$ is parallel. Then the sectional curvature $K$ of $M^n$ satisfies $K\geq 0$ if $c\geq 0$ and $K\geq c$ if $c<0.$
\end{proposition}

Combining Proposition \ref{prop} with a result of Erbacher  (\cite{Erb}, theorem $1$) we immediately obtain the following classification for type $(1, n-1)$ isometric immersions with constant scalar and mean curvatures (cf. \cite{DT}, theorem 1; also \cite{Sm}, \cite{CY}, \cite{Ros}). In contrast to previous authors we do not assume the non-negative sectional curvature condition in advance, but deduce it as a consequence of our other assumptions.
\vspace{2mm}

 \begin{corollary}\label{cor}
 	Let $M^n, n\geq 3,$ be a complete Riemannian manifold with constant scalar curvature and $f: M^n\rightarrow M^{n+k}(c), c\geq 0,$ be an isometric immersion with flat normal bundle which is proper and has type $(1, n-1).$ Assume that the mean curvature field of $f$ is parallel. Then $f(M^n)$ must be of the form ${\mathbb S}^1(r_1)\times {\mathbb S}^{n-1}(r_2),\  {\mathbb S}^1(r)\times {\mathbb R}^{n-1},$ or $\mathbb R\times {\mathbb S}^{n-1}(r)$ if $c=0;$ and $f(M^n)$ must be of the form ${\mathbb S}^1(r_1)\times {\mathbb S}^{n-1}(r_2)$ if $c>0.$ The corresponding local result holds if the completeness assumption on $M^n$ is dropped. 
 	\end{corollary}

 \begin{remark}\label{remark1}
 	If $f$ has parallel mean curvature field then the flat normal bundle assumption in Proposition \ref{prop} is redundant when $k=2$ and $M^n$ is not minimal (\cite{Erb}, theorem $1'$).
 \end{remark}
 
 %%%%%%%%%%%%%%%%%%%%%%%%%%%%%%%%%%%%%%%%%%%%%%%%%%%%%%%%%%%%%%%%%%
 %%%%%%%%%%%%%%%%%%%%%
In $\S{5}$ we specialize to hypersurfaces and derive an analogue of Theorem \ref{theorem} for the higher order mean curvatures of a hypersurface. Let $f: M^n\rightarrow M^{n+1}(c)$ be an orientable hypersurface and $r, 1\leq r \leq n,$ be any integer. The $r$-th higher order mean curvature $H_r$ of $f$ is defined by $$^n C_r H_r=\sum_{1\leq i_1\cdots <i_r\leq n} \lambda_{i_1}\cdots \lambda_{i_r},$$ where $\{\lambda_1, \cdots, \lambda_n\}$ denotes the set of principal curvatures of $f$ with respect to a chosen smooth unit normal field on $M^n.$ In particular, $H_1=\frac{\lambda_1+\cdots +\lambda_n}{n}$ is the mean curvature of $f$ and $H_n=\lambda_1\cdots\lambda_n$ is the Gauss-Kronecker curvature of $f.$ It is well-known that $H_r$ is intrinsic when $r$ is even, and that $|H_n|$ is intrinsic for any $n.$ Hypersurfaces with constant higher order mean curvatures are extensively studied (cf. \cite {Rei}, \cite{Rose}, \cite{AM}, \cite{NZ}, \cite{MP}, \cite{Ros}, \cite{Ros1}).
%(Of course, $H_1$ is not intrinsic; are $H_r, r\geq 3,$ intrinsic? Possible (See petersen p. 99, Allendorfer rigidity theorem in Kob-Nom and spival vol 5). 
\vspace{2mm}

The proof of Theorem \ref{theorem} can be readily adopted to obtain the following result which extends a result of Wu (\cite{Wu}, theorem 5.2, part (2)) to $c\neq 0$ when $r<n.$
\vspace{2mm}

\begin{theorem}\label{p1}
	Assume that $f: M^n\rightarrow M^{n+1}(c), n\geq 4,$ is a complete conformally flat orientable hypersurface with constant $r$-th higher order mean curvature $H_r.$ If $r$ is even and $r<n,$ then $H_r\geq 0.$ 
\end{theorem}
\vspace{2mm}

	In Theorem \ref{p1}, suppose $M^n$ is compact. If $c=0,$ there exists a point $x\in M^n$ such that the principal curvatures of $M^n$ at $x$ are all strictly positive with respect to a suitably chosen normal. Consequently, $H_r>0$ for each even $r.$ Same conclusion holds for $c<0,$ but clearly not for $c>0.$
	
	%as well. If $c=-1,$ there exists a point $p\in M$ where all the principle curvatures are greater than  $1$ (\cite{A}, lemma 8), and in particular $H_r>0$ for each even $r.$

\vspace{2mm}

Finally, we show that our earlier result (\cite{G}, theorem 1.1) extends to conformally flat hypersurfaces with constant $r$-th higher order mean curvature $H_r$ for any $r.$

\begin{proposition}\label{p2}
	Let $f: M^n\rightarrow M^{n+1}(c), n\geq 4,$ be a connected conformally flat orientable hypersurface of a space form $M^{n+1}(c).$ Assume that $M^n$ has constant $r$-th higher order mean curvature $H_r$ for some integer $r, 1\leq r\leq n.$  Then the following holds.
	
	\noindent (i) Suppose $H_r\neq 0.$ Then either $M^n$ is totally-umbilical or it has no umbilical points and $f(M^n)$ is an open subset of a rotation hypersurface in $M^{n+1}(c).$\\
		\noindent (ii) Suppose $H_r=0.$ Then either $M^n$ has constant sectional curvature $c$ or it has no umbilical points and $f(M^n)$ is an open subset of a rotation hypersurface in $M^{n+1}(c).$ 
\end{proposition}

\vspace{2mm}
 
%%%%%%%%%%%%%%%%%%%%%%%%%%%%5
%%%%%%%%%%%%%%%%%%%%%%%%%%%%%%
\section{Preliminaries}
\bigskip

Assume that $M^n$ and $f$ are as in Theorem \ref{theorem}. Let $\eta_1(x)$ and $\eta_2 (x)$ denote the principal normal vectors at $x$ in $M^n$ with multiplicities $1$ and $n-1,$ respectively. In this case, it is well-known that the distributions $E_{\eta_1}$ and $E_{\eta_2}$ are integrable and the immersion $f$ is locally holonomic. Thus each point $p$ in $M^n$ admits a principal coordinate system $(U, (x_1, \cdots, x_n))$ such that $E_{\eta_1}=span \{\partial_1\}$ and $E_{\eta_2}=span \{\partial_2, \cdots,  \partial_n\}$ on $U.$ In these coordinates, the first and the second fundamental forms of $f$ are given respectively by 

$$ds^2=v_1^2dx_1^2+\cdots+v_n^2dx_n^2$$ 
 \noindent and

$$\alpha(x) (\partial_1, \partial_1)=v_1^2\eta_1 (x)$$

$$\alpha(x) (\partial_1, \partial_j)=0$$

$$ \alpha(x) (\partial_i, \partial_j)=\delta_{ij}v_j^2\eta_2 (x),$$
\smallskip

 \noindent where $2\leq i,  j\leq n$ and $v_i=v_i(x_1,\cdots, x_n),  1\leq i\leq n,$ are some positive functions defined on $U.$
 \vspace{2mm}
 %%%%%%%%%%%%%%%%%%
 
Applying the Gauss equation to the orthonormal set 
$$\{X_i=\frac{\partial_i}{v_i}, X_j=\frac{\partial_j}{v_j}\} \ \ \  (1\leq i\neq j\leq n)$$ 

\noindent and making use of the expression for the second fundamental form of $f$ given above yields
\begin{equation}\label{g1}
K(X_1, X_j)=c+\langle \eta_1, \eta_2\rangle\ \ \ (2\leq j\leq n)
\end{equation}
\noindent and 
\begin{equation}\label{g2}
K (X_i, X_j)=c+||\eta_2||^2\ \ \ (2\leq i\neq j\leq n).
\end{equation}

\noindent Here $K(X, Y)$ denotes the sectional curvature of $M^n$ along the $2$-plane spanned by the vectors $X$ and $Y.$\vspace{2mm}

%%%%%%%%%%%%%%%%%%%%%
%$$\langle X_j, Y_j\rangle \nabla^{\perp}_{X_i}\eta_j=\langle \nabla_{X_j} Y_j, X_i\rangle (\eta_j-\eta_i) \ \ \  \bigl(X_i\in \Gamma(E_{\eta_i}),  X_j, Y_j\in \Gamma(E_{\eta_j}) \ \ (1\leq i\neq j\leq n)\bigr).$$

For any vector bundle $E,$ let $\Gamma(E)$ denote the set of all smooth local sections of $E.$ In the current setting, the Codazzi equations take the form (\cite{DT1}, p. 36, equation (1.41))
$${\nabla}^{\perp}_{X} \eta_2=0 \ \ \  (X\in \Gamma(E_{\eta_2}))$$

$$\langle X, Y\rangle \nabla^{\perp}_{Z}\eta_2=\langle \nabla_{X} Y, Z\rangle (\eta_2-\eta_1) \ \ \  \bigl(Z\in \Gamma(E_{\eta_1});  X, Y\in \Gamma(E_{\eta_2})\bigr)$$

\noindent and 
$$\langle X, Y\rangle \nabla^{\perp}_{Z}\eta_1=\langle \nabla_{X} Y, Z\rangle (\eta_1-\eta_2) \ \ \  \bigl(Z\in \Gamma(E_{\eta_2});  X, Y\in \Gamma(E_{\eta_1})\bigr).$$\vspace{2mm}

\noindent Substituting coordinate vector fields for $X, Y$ and $Z$ in these equations and using the expression for the first fundamental form of $f$ given above we obtain
\begin{equation}\label{c1}
{\nabla}^{\perp}_{\partial_j} \eta_2=0
\end{equation}

\begin{equation}\label{c2}
{\nabla}^{\perp}_{\partial_1} \eta_2=\frac{\partial_1 v_j}{v_j}(\eta_1-\eta_2)
\end{equation}

\noindent and

\begin{equation}\label{c3}
{\nabla}^{\perp}_{\partial_j} \eta_1=\frac{\partial_j v_1}{v_1}(\eta_2-\eta_1),
\end{equation}

\noindent where $ 2\leq j\leq n.$
\vspace{2mm}
%%%%%%%%%%%%%%%%%%%%%%%

Let $D=\frac{R}{2(n-1)}.$ Since $M^n$ has scalar curvature $R$ we obtain 
\begin{equation}\label{csc}
c+\langle\eta_1, \eta_2\rangle +\frac{n-2}{2}\bigl(c+||\eta_2||\bigr)^2=D
\end{equation}
\noindent from (\ref{g1}), (\ref{g2}) and the Gauss equation.
\vspace{2mm}

%For later use in the proof of Theorem \ref{theorem}, we note that in the special case when $\partial_j v_1=0$ on $U$ for each $j, \ 2\leq j\leq n,$ we may assume that $v_1=v_1(x_1).$ After fixing some point $P_0\in U$ and making the change of variable 
%$$t=t(x_1)=\int_{x_1(P_0)}^{x_1} v_1(u)du,$$

%\noindent we may assume that the first fundamental form of $f$ in the resulting new principal coordinates $(t, x_2, \cdots, x_n)$ is given by
%$$ds^2=dt^2+v_2^2dx_{2}^2+\cdots+ v_n^2 dx_{n}^2.$$
%\vspace{2mm}

% Finally, we assume that $M^n$ is $C^{\infty}$ and  connected.

%%%%%%%%%%%%%%%%%%%%%%%%%%%%%%%%%%%%%%%%%%%%%%%%%%%%%%%%%%%%%%%%%%%%%%%%%%%%%%%%%%%%%%

\section{A Scalar Curvature Inequality}
\medskip

This section is devoted to the proof of Theorem \ref{theorem}. We begin with a few preparatory lemmas.
\vspace{2mm}

Define two functions $\psi, \phi: M^n\rightarrow \mathbb R$ by
$$\phi(x)=c+\langle \eta_1 (x), \eta_2(x)\rangle\ \ \  (x\in M^n)$$

\noindent and 
$$\psi=\frac{2D}{n}-\phi.$$
\smallskip

\noindent Let $(U, (x_1, \cdots, x_n))$ be any connected principal coordinate chart in $M^n.$ Since $D$ is constant by assumption, using (\ref{c2}) and (\ref{csc}) we compute \vspace{2mm}

\vspace{-\baselineskip}
\begin{align*}
\partial_1 \phi=&\partial_1 \bigl(D-\frac{n-2}{2}(c+||\eta_2||^2)\bigr)\\
&= -(n-2) \langle \nabla_{\partial_1}^{\perp} \eta_2, \eta_2\rangle\\
&=-(n-2) \langle \frac{\partial_1 v_j}{v_j}(\eta_1-\eta_2), \eta_2\rangle\\
&=-(n-2)\frac{\partial_1 v_j}{v_j}\bigl((c+\langle \eta_1, \eta_2\rangle)-(c+|\eta_2|^2)\bigr)\\
&=-(n-2)\frac{\partial_1 v_j}{v_j}\bigl(\phi-(c+|\eta_2|^2)\bigr)\\
&=n\frac{\partial_1 v_j}{v_j}(\frac{2D}{n}-\phi).
\vspace{-\baselineskip}
\end{align*}

\noindent Therefore
\begin{equation}\label{ode}
\partial_1 \psi+n \frac{\partial_1 v_j}{v_j} \psi=0.
\end{equation}
\vspace{2mm}

\noindent Similarly, from (\ref{c1}) it follows that 
\begin{equation}\label{jpartials}
\partial_j \phi=0, \ 2\leq j\leq n.
\end{equation}\vspace{2mm}

\noindent On the other hand, for any $j$ with $2\leq j\leq n,$ using (\ref{c1}) and (\ref{c3}) we obtain \vspace{2mm}

\vspace{-\baselineskip}

\begin{align*}
	0&= \partial_j \phi\\
	&=\langle \nabla^{\perp}_{\partial_j}\eta_1, \eta_2\rangle\\
	& = \langle \frac{\partial_j v_1}{v_1}(\eta_2-\eta_1), \eta_2 \rangle.
\vspace{-\baselineskip}
\end{align*}

\noindent Thus 
\begin{equation}\label{jpartials1}
\frac{\partial_j v_1}{v_1} ({||\eta_2||}^2-\langle \eta_1, \eta_2\rangle)=0\ \ on \ \ U.
\end{equation}\vspace{2mm}

Recall from $\S{2}$ that $D:=\frac{R}{2(n-1)}$ and $X_i:=\frac{\partial_i}{v_i}\  (1\leq i\leq n).$\vspace{2mm}

\begin{lemma}\label{L1}
	If ${||\eta_2||}^2(p)=\langle \eta_1, \eta_2\rangle (p)$ for some $p\in M^n,$ then $M^n$ has constant sectional curvature and $D=\frac{nc}{2}.$
\end{lemma}

\noindent {\bf Proof.} Since $\bigl({||\eta_2||}^2-\langle \eta_1, \eta_2\rangle \bigr)(p)=0,$ we have $\phi(p)=c+|\eta_2(p)|^2.$ Therefore (\ref{csc}) gives $\psi(p)=0.$
From (\ref{ode}) and (\ref{jpartials}) we get $\psi(x)=0$ for every $x\in U,$ and consequently $\phi(x)=\frac{2D}{n}$ for every $x\in M^n$ by covering $M^n$ by a collection of overlapping principal coordinate charts. Since $f$ has flat normal bundle and the curvature operator of $M^n$ is diagonalized by the orthonormal basis $\{X_i\wedge X_j: 1\leq i < j\leq n\},$ it follows that $M^n$ has constant sectional curvature $\frac{2D}{n}$ (\cite{P}, p. 63). Further, since $\frac{2D}{n}=c+||\eta_2||^2$ we conclude $D= \frac{nc}{2}$ in this case. \qed

\vspace{2mm}

%%%%%%%%%%%%%%%%%%%%%%%%%%%%%%%%%%%%%%%%
In view of Lemma \ref{L1} we may assume that $\bigl({||\eta_2||}^2-\langle \eta_1, \eta_2\rangle\bigr) (p)\neq 0$ for every $p\in M^n.$ From (\ref{jpartials1}) $\partial_jv_1=0$ for each $j, 2\leq j\leq n.$ After fixing a point $p_0\in U$ and making the change of variable 
 $$t=t(x_1):=\int_{x_1(p_0)}^{x_1} v_1(u)du,$$
 
 \noindent we may assume that the first fundamental form of $f$ in the resulting principal coordinates $(t, x_2, \cdots, x_n)$ is given by
 $$ds^2=dt^2+v_2^2dx_{2}^2+\cdots+ v_n^2 dx_{n}^2.$$ %As noted in $\S{2},$ it  follows that each point of $M^n$ admits a principal coordinate system $(U, (t, x_2, \cdots, x_n))$ in which the first fundamental form of the immersion $f$ is given by
%$$I=dt^2+v_2^2dx_{2}^2+\cdots+ v_n^2 dx_{n}^2.$$

\noindent In the remainder of this proof we assume that this reduction has already been made and only work with principal coordinate systems on $M^n$ that are of this form.\vspace{2mm}

%%%%%%%%%%%%%%%%%%%%%%%%%%%%%%%%%%%%%%%
\begin{lemma}\label{L2} 
	We have $K(\partial_t, \partial_j)=-\frac{1}{v_j} \partial_{tt}v_j$ for any $j$ with $2\leq j\leq n.$
	\end{lemma}
	
 \noindent {\bf Proof.} This is standard. See, for example, (\cite{DD}, p. 20).
 \qed
 %%%%%%%%%%%%%%%%%%%%%%%
 %%%%%%%%%%%%%%%%%%%%%%%%%%
 % but for the sake of completeness we outline the proof here. Since $||\partial_t||=1$ we have $\nabla_{\partial_t} \partial_t=0.$  We set $X_j:=\frac{\partial_j}{v_j}, h_{ji}:=\frac{\partial_j v_i}{v_j}$ and $h_{1j}:=\partial_t v_j\  (2\leq i\neq j\leq n).$ A straightforward computation (see \cite{DT1}, p. 20) gives 
 
%$$\nabla_{\partial_j} \partial_t=h_{1j}X_j$$

%\noindent and 

%$$\nabla_{\partial_t} X_j=0$$ 

%\noindent for any $j, 2\leq j\leq n.$
%\vspace{2mm}

%Let $\mathcal R$ denote the Riemann curvature tensor on $M^n.$ Using the expressions of the preceeding paragraph we compute 

%\vspace{-\baselineskip}
%\begin{align*}
%&\mathcal R(\partial_j, \partial_t)\partial_t\\
%&=\nabla_{\partial j} \nabla_{\partial_t} \partial_t-\nabla_{\partial_t}\nabla_{\partial j} \partial_t\\
%&=-\nabla_{\partial_t}\nabla_{\partial j} \partial_t\\
%&=-\nabla_{\partial_t}(h_{1j}X_j)\\
%&=-\partial_t (h_{1j}) X_j-h_{1j}\nabla_{\partial t} X_j\\
%&=-\partial_t (h_{1j}) X_j.
%\vspace{-\baselineskip}
%\end{align*}

%\noindent Therefore 
%$$K(\partial_t, \partial_j)=-\frac{1}{v_j} \partial_{tt}v_j.$$
%%%%%%%%%%%%%%%%%%%%%%%%%%%%%%%%%%%%%%%%%
%%%%%%%%%%%%%%%%%%%%%%%%%%%%%%%%%%%%%%%%%%

\vspace{2mm}
%%%%%%%%%%%%%%%%%%%%%%%%%%%%%%%%%%%%%%%%%%%%%%%%%%%%%

Fix a point $p\in M^n$ and a principal coordinate system $(U, (t, x_2,\cdots, x_n))$ around $p.$ Let $\gamma:\mathbb R\rightarrow M^n, \gamma=\gamma(s),$ denote the geodesic with initial data $\gamma (0)=p$ and $\gamma'(0)=\partial_t (p).$ Note that the geodesic is defined over $\mathbb R$ since $M^n$ is geodesically complete. Since the integral curves of $\partial_t$ are geodesics parametrized by arc-length, it follows from uniqueness that $${\gamma}'(s)=\partial_t (\gamma(s))$$ whenever $\gamma(s)\in U.$ If $(V, (\tau, y_2,\cdots, y_n))$ denotes another principal coordinate system that also contains $\gamma(s_0)$ for some $s_0\in \mathbb R$ then, since $\partial _{\tau}\in E_{\eta_1(\gamma(s_0))}$ and $E_{\eta_1(\gamma(s_0))}$ is $1$-dimensional, we must have $\partial_{\tau} (\gamma(s_0))=\pm \partial_{t} (\gamma(s_0)).$ Replacing $\tau$ by $-\tau$ if necessary, we can assume that $\partial_{\tau} (\gamma(s_0))=\partial_{t} (\gamma(s_0)).$ In this way we obtain a covering of ${\gamma}(\mathbb R)$ by principal coordinate charts $\{(U, (t, x_2,\cdots, x_n))\}$ which satisfy the condition $${\gamma}'(s)=\partial_t (\gamma(s))$$ whenever $\gamma(s)\in U.$\vspace{2mm}
%%%%%%%%%%%%%%%%%%%%%%%%%%
%%%%%%%%%%%%

Assuming that the function $\psi$ is positive for the moment, we define a function $\theta:\mathbb R\rightarrow \mathbb R$ by
 $$\theta(s)=\frac{1}{\psi(\gamma(s))^{\frac{1}{n}}}.$$ \vspace{2mm}

 \begin{lemma}\label{L3}
 	The function $\theta$ satisfies the equation $\frac{\theta''(s)}{\theta(s)}=-(\phi\circ \gamma)(s)$ for each $s\in \mathbb R.$
\end{lemma}

\noindent {\bf Proof.} Working in a principal coordinate system we obtain 
$$0=(\log(\psi\circ\gamma){\theta}^n)'= \frac{\partial_t \psi}{\psi}+n \frac{\theta'}{\theta}.$$

\noindent Here $'$ denotes derivative with respect to the arc-length parameter $s$ of $\gamma.$ Using 
$$\partial_t \bigl( \frac{\partial_t v_j}{v_j}\bigr)=\frac{1}{v_j}\partial_{tt} v_j-(\frac{\partial_t v_j}{v_j})^2$$

\noindent and Lemma \ref{L2} we obtain
$$\partial_t \bigl( \frac{\partial_t v_j}{v_j}\bigr)=-K(\partial_t, \partial_j)-(\frac{\partial_t v_j}{v_j})^2=-\phi\circ \gamma-(\frac{\partial_t v_j}{v_j})^2$$

\noindent and

$$\bigl(\frac{\theta'}{\theta}\bigr)'=-\phi\circ \gamma-\bigl(\frac{\partial_t v_j}{v_j}\bigr)^2=-\phi\circ \gamma-\bigl(\frac{\theta'}{\theta}\bigr)^2.$$

\noindent Therefore 
$$\frac{\theta''}{\theta}=-\phi\circ \gamma.$$  \qed

\vspace{2mm}

%%%%%%%%%%%%%%%%%%%%%
%%%%%%%%%%%%%%%%%
Now we consider the three cases separately and complete the proof of Theorem \ref{theorem}.\vspace{2mm}

\noindent (i) Let $c=0.$ We assume that $D<0$ and obtain a contradiction. Using (\ref{csc}) we obtain
$$\psi=\frac{2D}{n}-\phi=\frac{n-2}{2}||\eta_2||^2-\frac{n-2}{n}D\geq -\frac{n-2}{n}D>0.$$ 

\noindent Also, since $\phi<D,$ from Lemma \ref{L3} it follows that $\frac{\theta''}{\theta}> -D>0.$ %from which we conclude $$(e^{Ds}\theta)'(s)>0$$ and $$\theta(s)> a e^{-Ds}$$ for some positive constant $a.$ 
Since $\theta$ is strictly convex on $\mathbb R$ we have $\lim_{s\rightarrow \infty} \theta(s)=\infty.$ This contradicts $\psi(\gamma(s)){(\theta(s))}^n=1$ for all $s\in \mathbb R.$ Hence $D\geq 0.$ \vspace{2mm}
 
 Now assume that $M^n$ is compact and $D=0.$ Using (\ref{csc}) we obtain $\psi\geq 0.$ Assume that $\psi$ attains minimum at $q\in M^n.$ We claim that $\psi(q)=0.$ If $\psi(q)>0,$ then from (\ref{ode}) we obtain $\partial_t v_j(q)=0.$ The equation 
 
 $$\partial_{tt} \psi=-n\partial_t\bigl(\frac{\partial_t v_j}{v_j}\bigr)\psi+n^2\bigl(\frac{\partial_t v_j}{v_j}\bigr)^2,$$ 
 
 \noindent which is obtained from (\ref{ode}) by differentiation, and the condition $\partial_{tt} \psi(q)\geq 0$ together imply that $\partial_t\bigl(\frac{\partial_t v_j}{v_j}\bigr)(q)\leq 0.$ Using $\partial_t v_j(q)=0,$ (\ref{g1}) and Lemma \ref{L2} we conclude
  $$\phi(q)=K(\partial_t, \partial_j)(q)=-\frac{1}{v_j} \partial_{tt}v_j (q)= -\partial_t\bigl(\frac{\partial_t v_j}{v_j}\bigr)(q)\geq 0.$$
  
\noindent Hence $\psi(q)=-\phi(q)\leq 0$ and $\psi(q)=0$ which contradicts the assumption. Therefore $\psi(q)=0.$ From (\ref{ode}) and $\partial_j \phi=0$ we conclude that $\phi$ vanishes on $M^n.$ From (\ref{csc}) we get $||\eta_2||^2=0$ and it follows that $M^n$ must be flat. Since $\eta_2$ vanishes identically, $f$ has constant index of relative nullity $n-1.$ From Hartman's theorem \cite{H} $f$ must be an $(n-1)$-cylinder and $M^n$ contains a line. This contradicts the compactness of $M^n.$\vspace{2mm}

\noindent (ii) Let $c>0.$ Suppose $D\leq \frac{n-2}{2}c.$ From ({\ref{csc}) we have $$\psi=\frac{2D}{n}-\phi=\frac{n-2}{2}c+\frac{n-2}{2}||\eta_2||^2\geq \frac{n-2}{2}c>0.$$
	%=\frac{n-2}{n}D\geq \frac{n-2}{2n}(nc-2D)>\frac{n-2}{2n}c>0.$$ 
\noindent Moreover, (\ref{csc}) implies that $\phi+\frac{n-2}{2}c\leq D,$ or $\phi\leq D-(\frac{n-2}{2}c)\leq 0.$ From Lemma \ref{L2} it follows that $\theta$ is convex, and strictly convex if $D<\frac{n-2}{2}c.$ As in part (i) we obtain a contradiction if $D<\frac{n-2}{2}c,$ or if $D=\frac{n-2}{2}c$ and one of the limits $\lim_{s\rightarrow \pm \infty} \theta (s)$ is infinite.  Otherwise $\theta$ must be a constant function and $\psi,$ and hence $\phi,$ must be constant along $\gamma.$ Since $\partial_j \phi=0, \ 2\leq j\leq n,$ it follows that $\phi$ must be a constant function on $M^n.$ From (\ref{ode}) we obtain $\partial_t v_j=0$ for each $j,\ 2\leq j\leq n.$ Hence  $\phi=K(\partial_t, \partial_j)=0$ using Lemma \ref{L2}. Since $D=\frac{n-2}{2}c,$ using (\ref{csc}) we obtain $\eta_2=0.$ This gives $0=\phi=c+\langle \eta_1, \eta_2\rangle=c,$ which is a contradiction.
\vspace{2mm}

\noindent (iii) Let $c<0.$ Suppose $D<\frac{nc}{2}.$ Since $nc-2D>0,$ from (\ref{csc}) we obtain $$\psi =\frac{2D}{n}-\phi>\frac{n-2}{2n}(nc-2D)>0.$$ Using (\ref{csc}) once more we get $\phi+\frac{n-2}{2}c\leq D,$ or $\phi\leq D-(\frac{n-2}{2}c)<c<0$ since $D<\frac{nc}{2}$ by assumption. Rest of the proof is analogous to part (i).

In the case of equality, we have $\phi+\frac{n-2}{2}{||\eta_2||}^2=c<0$ so that $\phi(x)<0$ for every $x\in M^n.$ However, if $M^n$ is compact then the argument in part (i) shows that there exists a point $q\in M^n$ such that $\phi(q)\geq 0.$ Hence the inequality must be strict if $M^n$ is compact.
\qed
\vspace{2mm}
%%%%%%%%%%%%%%%%
%%

%%%%%%%%%%%%%%%%%%%%%%%%%%%%%%%%%%%%
%%%%%%%%%%%%%%%%%%%%%%%%%%%%%%%%%%%%%
\section{A Sectional Curvature Inequality}
\medskip

In this section we prove Proposition \ref{prop}.\vspace{2mm}

 As in the proof of Theorem \ref{theorem} we work in a principal coordinate system $(U, (x_1, \cdots, x_n))$ on $M^n.$ We claim that $\partial_1 v_j$ vanishes identically for every $j, 2\leq j\leq n.$  Otherwise, by shrinking $U$ if necessary, we may assume that  $\partial_1 v_j $ does not vanish on $U.$  Let $H$ be the mean curvature field of $f.$ Since $nH=\eta_1+(n-1)\eta_2$ and $\nabla^{\perp} H=0,$ from (\ref{c1}) and (\ref{c3}) it follows that  
\begin{equation}\label{m1}
\nabla_{\partial_j}^{\perp} \eta_1=0
\end{equation}

\noindent and 
\begin{equation}\label{m2}
  \nabla^{\perp}_{\partial_1} \eta_1+(n-1)
\nabla^{\perp}_{\partial_1} \eta_2=0.
\end{equation}

 \noindent Thus $\partial_jv_1=0$  for  each $j, 2\leq j\leq n.$ From (\ref{c2}) we conclude that
\begin{equation}\label{m3}
 \frac{\partial_1v_j}{v_j}=\frac{\partial_1v_k}{v_k}\ \ \  (2\leq j, k\leq n).
 \end{equation}
  
\noindent Since $M^n$ has flat normal bundle, using (\ref{c1}) we obtain

\vspace{-\baselineskip}
%\begin{align}

%$0={\mathcal R}^{\perp}(\partial_1, \partial_j)\eta_2$

%$=\nabla^{\perp}_{\partial_1} (\nabla^{\perp}_{\partial_j}\eta_2)-\nabla^{\perp}_{\partial_j} (\nabla^{\perp}_{\partial_1}\eta_2)$

%$=-\nabla^{\perp}_{\partial_j} (\nabla^{\perp}_{\partial_1}\eta_2).$
%\vspace{-\baselineskip}
%\end{align}
%%%%%%%%%%%%%%%%%%%%%%%%%%%%%%%%%%%%%%%%%%%%%%%%%%%%%%%%
\begin{align*}
0& ={\mathcal R}^{\perp}(\partial_1, \partial_j)\eta_2 \\
&=\nabla^{\perp}_{\partial_1} (\nabla^{\perp}_{\partial_j}\eta_2)-\nabla^{\perp}_{\partial_j} (\nabla^{\perp}_{\partial_1}\eta_2) \\
&=-\nabla^{\perp}_{\partial_j} (\nabla^{\perp}_{\partial_1}\eta_2).
\end{align*}
%%%%%%%%%%%%%%%%%%%%%%%%%%%%%%%%%%%%%%%%%%%%%%%%%%%%%%%%%
\noindent Substituting for $\nabla^{\perp}_{\partial_1} \eta_2$ from (\ref{c2}) and simplifying we obtain
$$\nabla^{\perp}_{\partial_j}\eta_1=\frac{\partial_j(\frac{\partial_1 v_j} {v_j})}{\frac{\partial_1 v_j} {v_j}} (\eta_2-\eta_1).$$

\noindent From (\ref{m1}) it follows that 
$$\partial_j(\frac{\partial_1v_j}{v_j})=0 \ {for\ each\ } j, 2\leq j\leq n.$$
	
	\noindent Combining this with equation (\ref{m3}) we obtain $$\partial_k (\frac{\partial_1v_j}{v_j})=0$$ for any $k, 1\leq k\leq n.$ Hence we may assume that
	$$v_k(x_1, \cdots, x_n)=\psi_k (x_2, \cdots x_n) \mu(x_1)$$
	
	\noindent for some positive functions $\psi_k$ and $\mu.$ Differentiating (\ref{csc}) with respect to $x_1$ we obtain
	 $$\langle \nabla^{\perp}_{\partial_1} \eta_1, \eta_2\rangle+\langle \eta_1, \nabla^{\perp }_{\partial_1}\eta_2\rangle+(n-2) \langle \eta_2, \nabla^{\perp}_{\partial_1}\eta_2\rangle=0.$$
	
	\noindent Substituting $\nabla^{\perp}_{\partial_1}\eta_1=-(n-1)\nabla^{\perp}_{\partial_1}\eta_2$ from (\ref{m2}) in the above equation gives
	$$ \langle \nabla^{\perp}_{\partial_1}\eta_2, \eta_1-\eta_2\rangle=0.$$
	
\noindent 	Since $\nabla^{\perp}_{\partial_1}\eta_2=\mu'(x_1) (\eta_1-\eta_2)$ we have $\mu'(x_1) ||\eta_1-\eta_2||^2=0$ so that $\mu'$ vanishes identically and hence ${\partial_1} v_j=0.$ From this and $\partial_jv_1=0$ we conclude $K(\partial_1, \partial_j)=0$ (\cite{DD}, p. 20). Since the curvature operator of $M^n$ is diagonalized by the orthonormal basis $\{\frac{\partial_i}{v_i}\wedge \frac{\partial_j}{v_j}: 1\leq i<j\leq n \},$ we conclude as in Lemma \ref{L1} that the sectional curvatures of $M^n$ lie between $0$ and $c+||\eta_2||^2.$ 
	\qed

\bigskip

\section {Hypersurfaces and higher order mean curvature}
\medskip

  In this section we specialize to hypersurfaces and prove Theorem \ref{p1} and Proposition \ref{p2}.
  \vspace{2mm}
  
  Let $f: M^n\rightarrow M^{n+1}(c)$ be an orientable hypersurface and $r, 1\leq r \leq n,$ be any integer. The $r$-th higher order mean curvature $H_r$ of $f$ is defined by $$^n C_r H_r=\sum_{1\leq i_1\cdots <i_r\leq n} \lambda_{i_1}\cdots \lambda_{i_r},$$ where $\{\lambda_1, \cdots, \lambda_n\}$ denotes the set of principal curvatures of $f$ with respect to a chosen smooth unit normal field on $M^n.$ It is well-known that $H_r$ is intrinsic when $r$ is even, and that $|H_n|$ is intrinsic for any $n.$ 
	
%	(******* Are $|H_r|'s$ intrinsic ? $H_r$ is intrinsic when $r$ is even, and so is $H_n.$ Are $H_r, r\geq 3$ also intrinsic? Possible.}

\vspace{2mm}

%%%%%%%%%%%%%%%
%%%%%%%%%%%%%%%%%%
	
%	Result is optimal when $c=0.$ See O. Palmas, Rotation hypersurface with constant $H_r.$

\noindent \textit{Proof of Theorem \ref{p1}}. We only need to modify the proof of Theorem \ref{theorem} slightly. If $f$ is umbilical at some point $p\in M^n$ then the conclusion is trivial. Otherwise, since $M^n$ is conformally flat and $n\geq 4,$ by the aforementioned result of Cartan \cite{C} we may assume that the principal curvatures of $f$ satisfy the relation $\lambda_1= \cdots= \lambda_{n-1}:=\lambda$ and $\lambda_n:=\mu,$ where $\lambda$ and $\mu$ denote the \textit{distinct} principal curvatures of $M^n$ with respect to a smooth unit normal field $\mathcal{N}$ on $M^n.$ Let $\eta_1:=\mu {\mathcal N}$ and $\eta_2:= \lambda \mathcal {N}.$ Then $\eta_1$ and $\eta_2$ are the principal normal vector fields of $f$ and $f$ has type $(1, n-1).$ From here onwards we proceed as in the proof of Theorem \ref{theorem}. The Codazzi equations can be written in terms of the principal curvatures as
\begin{equation}\label{c15}
\partial_j (\lambda)=0
\end{equation}

\begin{equation}\label{c25}
\partial_1 (\lambda)=\frac{\partial_1 v_j}{v_j} (\mu-\lambda)
\end{equation}

\begin{equation}\label{c35}
\partial_j (\mu)=\frac{\partial_j v_1}{v_1}(\lambda_\mu)
\end{equation}

\noindent for every $j, 1\leq j\leq n.$
\vspace{2mm}

By definition
$$^nC_rH_r=\sum_{1\leq i_1\cdots <i_k\leq n} \lambda_{i_1}\cdots \lambda_{i_r} =^{n-1}C_r {\lambda}^r+^{n-1}C_{r-1}{\lambda}^{r-1} \mu.$$

\noindent Let 
\begin{equation}\label{D}
D:=\frac{n}{{r}}H_r=\frac{n-r}{r}{\lambda}^r+{\lambda}^{r-1} \mu.
\end{equation}
\vspace{2mm}

Define functions $\psi, \phi: M^n\rightarrow \mathbb R$ by
$$\phi={\lambda}^{r-1}\mu$$

\noindent and
$$\psi=\frac{r}{n}D-\phi.$$

\noindent Note that $\phi$ is well-defined since $r$ is even.
\vspace{2mm}

Let $(U, (x_1, \cdots, x_n))$ be any principal coordinate system in $M^n.$ Since $D$ is constant by assumption, from (\ref{c15}) and (\ref{D}) it follows that 
\begin{equation}\label{jpartials5}
\partial_j \phi=0, \ 2\leq j\leq n.
\end{equation}

\noindent On the other hand, for any $j\geq 2,$ using (\ref{jpartials5}) and (\ref{c35}) we obtain

\vspace{-\baselineskip}
\begin{align*}
	0&= \partial_j \phi\\
	&={\lambda}^{r-1}\partial_j(\mu)\\
	& = {\lambda}^{r-1} \frac{\partial_j v_1}{v_1} (\lambda-\mu).
	\vspace{-\baselineskip}
\end{align*}

\noindent Thus 
\begin{equation}\label{jpartials1 5}
\frac{\partial_j v_1}{v_1} ({\lambda}^r-\phi)=0\ \ on \ \ U.
\end{equation}

\vspace{2mm}

\noindent If there exists a point $p\in M^n$ such that ${\lambda}^r (p)=\phi(p),$ then from (\ref{D}) it follows that $D\geq 0$ since $r$ is even. Therefore we may assume that $\frac{\partial_j v_1}{v_1}=0$ for all $2\leq j\leq n.$

	%In fact $\phi=frac{rD}{n}={\lambda}^r\geq 0.$ Curvature is constant by Gauss equation.
	
%%%%%%%%%5
%%%%%%%%%
%using (\ref{D}) and (\ref{c25}) we compute %(************this uses $r<n?$)\vspace{2mm}

%\vspace{-\baselineskip}
%\begin{align*}
%	\partial_1 \phi=&\partial_1 \bigl(D-\frac{n-r}{r}{\lambda}^r\bigr)\\
%	&=-(n-r)\frac{\partial_1 v_j}{v_j} (\phi-{\lambda}^r)\\
%	&=\frac{1}{n}\frac{\partial_1 v_j}{v_j}(\frac{r}{n}D-\phi).
%\vspace{-\baselineskip}
%\end{align*}

%\noindent Therefore
%\begin{equation}\label{ode5}
%\partial_1 \psi+n \frac{\partial_1 v_j}{v_j} \psi=0.
%\end{equation}
%\vspace{2mm}
%%%%5
%%%%%%5

\vspace{2mm}

Following the proof of Theorem \ref{theorem}, we note Lemmas \ref{L2} and \ref{L3} continue to hold without any changes with the same definition for $\theta$ as before.\vspace{2mm}

 We assume $D<0$ and obtain a contradiction. Since $D<0,$ using (\ref{D}) we conclude
  \begin{align*}
 \psi&=\frac{rD}{n}-\phi\\
 &= \frac{rD}{n}-D+\frac{n-r}{r}{\lambda}^r\\
& > \frac{r-n}{n}D\\
 &> 0.
\end{align*}
This, combined with strict convexity of $\theta,$ gives a contradiction as in case (i) of Theorem \ref{theorem}.
\qed

%%%%%%%%%%%%%%%%%%%%5
%\noindent {\bf Possible improvement added. Claim: $H_n\geq 0$ if $n$ is even and $c=0.$ To prove this, we consider the product immersion $F=(f, id_{\mathbb R}): M^n\times \mathbb R\rightarrow {\mathbb R}^{n+1}\times {\mathbb R}={\mathbb R}^{n+2}$ given by $F(x, t)=(f(x), t)$ for any $(x, t)\in M^{n}\times \mathbb R.$ Since $\lambda_{n+1} (F)=0,$ we have $H_n(f)=H_n(F)\geq 0$ from the above. Question: Can something similar be done for $c<0$ and $c>0?$ But this needs conformal flatness of the product $M^n\times \mathbb R.$ Warning! $M^n\times \mathbb R$ is conformally flat if and only if $M^n$ has constant sectional curvature. But the formulation with only two distinct principal curvatures may still work. See the example in (\cite{DD}, p. 545) for conformal flatness of the warped products.}
%\qed
%%%%%%%%%%%%%%%%%%%%%%%%%%%%%%%%%%%%%%%%%%%%%%%%
\vspace{2mm}
%%%%%%%%%%%%%%%%%%%%%%%%%%%%%%%%%%%%%%%5
%%%%%%%%%%%%%%%%%%%%%%%%%%%%%%%%%%%%%%%%%

\noindent \textit{Proof of Proposition \ref{p2}}. We may assume that $\lambda_1= \cdots=\lambda_{n-1}:=\lambda$ and $\lambda_n:=\mu$ denote the (not necessarily distinct) principal curvatures of $M^n.$ From here, we proceed as in the original proof (\cite{G}, Theorem $1.1$). It is sufficient to replace the function $\psi$ appearing in that proof by the function $$\Psi:={\lambda}^r-\lambda^{r-1}\mu.$$

\noindent Suppose $\lambda(s)\neq 0$ for all $s.$ Setting $'=\frac{d}{ds}$ and using ${\lambda}'=\frac{1}{2\lambda}({\lambda}^2)',$ we compute
\begin{align*}
	({\lambda}^r)' =&r{\lambda}^{r-1} {\lambda}'\\
=&\frac{r}{2}{\lambda}^{r-2}({\lambda}^2)'\\
=&\frac{r}{2}{\lambda}^{r-2}\bigl(2\frac{x'}{x}(\lambda\mu-{\lambda}^2)\bigr)\\
=& r\frac{x'}{x} ({\lambda}^{r-1}\mu-{\lambda}^r)\\
=&-r\frac{x'}{x}\Psi.
\end{align*}

\noindent In the third line of the previous computation, we have used the formulas (\cite{DD}, proposition 3.2)  $${\lambda}^2=\frac{\delta-cx^2-x'^2}{x^2}$$

\noindent and $$\mu=\frac{x''+cx}{\sqrt{\delta-cx^2-x'^2}}$$

\noindent to compute $({\lambda}^2)'.$
\vspace{2mm}

Assume $r<n.$ The definition of $H_r$ gives
$$^nC_r H_r=^{n-1}C_r{\lambda}^r+^{n-1}C_{r-1}{\lambda}^{r-1}\mu$$

\noindent Let $\overline{H_r}:=\frac{n}{n-r}H_r$ and $a:=\frac{r}{n-r}.$  Then 
$$\overline{H_r}={\lambda}^r+a{\lambda}.$$

\noindent Since $\overline{H_r}$ is constant using the above we obtain 

\begin{align*}
{\Psi}' =&({\lambda}^r-{\lambda}^{r-1}\mu)'\\
=& ({\lambda}^r+\frac{1}{a}{\lambda}^r-\frac{1}{a} \overline{H_r})'\\
=& (1+\frac{1}{a}) ({\lambda}^r)'\\
=& - (1+\frac{1}{a}) r \frac{x'}{x} \Psi.
\end{align*}

\noindent Noting $1+\frac{1}{a}=\frac{n}{r},$ we conclude 
$$\Psi'+n\frac{x'}{x}\Psi=0.$$ 

\noindent It follows that $$x^n\Psi=C$$ for some constant $C.$ We note that this relation is also valid for $r=n.$ Rest of the argument is the same as in \cite{G}.

\qed

\bigskip

%%%%%%%%%%%%%%%%%%%%%%%%%%%%
%%%%%%%%%%%%%%%%%%%%%%%%%%
 
 %%%%%%%%%%%%%%%%%%
 %%%%%%%%%%%%%%%%%%%%%%

\end{document}